\documentclass{amsart}
\usepackage{amscd}
\newtheorem{thm}{Theorem}[section]
\newtheorem{lemma}[thm]{Lemma}
\newtheorem{prop}[thm]{Proposition}
\newtheorem{cor}[thm]{Corollary}
\newtheorem{fact}[thm]{Fact}
\newtheorem{conj}[thm]{Conjecture}
\theoremstyle{definition}
\newtheorem{defn}[thm]{Definition}
\newtheorem{exam}[thm]{Example}
\newtheorem{prob}[thm]{Problem}
\newtheorem{quest}[thm]{Question}
\newtheorem{discuss}[thm]{Discussion}

\theoremstyle{remark}
\newtheorem{remark}[thm]{Remark}
\numberwithin{equation}{section}

\DeclareMathOperator{\Supp}{Supp}

\DeclareMathOperator{\Hom}{Hom}

\DeclareMathOperator{\Image}{Im}

\DeclareMathOperator{\codim}{codim}
\DeclareMathOperator{\indeg}{indeg}

\DeclareMathOperator{\reg}{reg}
\DeclareMathOperator{\Soc}{Soc}

\DeclareMathOperator{\emb}{embdim}
\DeclareMathOperator{\depth}{depth}
\DeclareMathOperator{\gr}{gr}
\DeclareMathOperator{\link}{link}

\newcommand{\frm}{{\mathfrak m}}
\newcommand{\fraM}{{\mathfrak M}}
\newcommand{\fraN}{{\mathfrak N}}
\newcommand{\frn}{{\mathfrak n}}
\newcommand{\frp}{{\mathfrak p}}

\newcommand{\Z}{\ensuremath{\mathbb Z}}

\newcommand{\bbR}{\ensuremath{\mathbb R}}
\newcommand{\bbZ}{\ensuremath{\mathbb Z}}
\newcommand{\bbinom}[2]{%
\genfrac{(}{)}{0pt}{}{#1}{#2}}

\newcommand{\sbinom}[2]{%
\text{\small $\displaystyle{\left(\!\!\!\begin{array}{c}#1 \\ #2 \end{array}\!\!\!
\right)}$}}

\begin{document}
\title{Buchsbaum homogeneous algebras with minimal multiplicity}
\author{Shiro Goto}
\address{Department of Mathematics, School of Science and Technology, 
Meiji University, 214--8571, Japan} 
\email{goto@math.meiji.ac.jp}
\curraddr{}
\author{Ken-ichi Yoshida}
\address{Graduate School of Mathematics, Nagoya University,
         Nagoya  464--8602, Japan}
\email{yoshida@math.nagoya-u.ac.jp}
\date{\today}
\dedicatory{}
\begin{abstract}
In this paper we first give a lower bound on multiplicities 
for Buchsbaum homogeneous $k$-algebras $A$ in terms of 
the dimension $d$, the codimension $c$, the initial degree $q$, and 
the length of the local cohomology modules of $A$. 
Next, we introduce the notion of Buchsbaum $k$-algebras 
with minimal multiplicity of degree $q$, and 
give several characterizations for those rings. 
In particular, we will show that those algebras have linear 
free resolutions. 
Further, we will give many examples of those algebras. 
\end{abstract}
\maketitle
\section{Introduction}

\par 
In this paper we are going to introduce a new class of 
non-Cohen--Macaulay Buchsbaum rings with linear resolution.  
Now let us explain our motivation. 
\par \vspace{2mm}
In 1967, Abhyankar \cite{Ab} proved that for a homogeneous integral domain $A$ 
over an algebraically closed field $k$, the following inequality  holds:
\[
 \emb(A)= \dim_k A_1 \le e(A) + \dim A -1,
\]
where $e(A)$ (resp. $\emb(A)$) denotes 
the multiplicity (resp. the embedding dimension) of $A$.  
In 1970, Sally \cite{Sa1} proved that the same inequality holds 
for any Cohen--Macaulay local ring, and called the ring $A$ 
which satisfies the equality 
a \textit{Cohen--Macaulay local ring with maximal embedding dimension}. 
In 1982, the first author \cite{Go2} pointed out that 
Sally's result can be extended to the class of Buchsbaum local rings. 
Namely, for any Buchsbaum local ring $A$, the following inequality holds:
\[
  \emb(A) \le e(A) + \dim A -1 +I(A),  
\]
where $I(A)$ denotes the Buchsbaum invariant ($I$-invariant) of $A$. 
A Buchsbaum local ring $A$ which satisfies the equality  is said to be 
a \textit{Buchsbaum local ring with maximal embedding dimension}. 
This is equivalent to that $A$ has $2$-linear resolution. 
\par 
On the other hand, the first author \cite{Go3} in 1983 proved that 
\[
 e(A) \ge 1 + \sum_{i=1}^{d-1} \bbinom{d-1}{i-1} l_A(H_{\frm}^i(A)),
\]
and defined the notion of 
\textit{Buchsbaum local rings with minimal multiplicity}.  
He also proved that those rings have $2$-linear resolutions. 

\par \vspace{2mm}
Recently, Terai and the second author \cite{TeYo1} defined the notion of 
\textit{Buchsbaum rings with minimal multiplicity with initial degree $q$
for Stanley--Reisner rings}, and proved that those rings have 
$q$-linear resolutions. 
Their result indicates the existence of Buchsbaum homogeneous $k$-algebras
with minimal multiplicity of higher initial degree. 

\par \vspace{2mm}
The main purpose of this paper is to introduce the notion of 
Buchsbaum homogeneous $k$-algebra with minimal multiplicity of
higher (initial) degree and characterize those rings. 
Moreover, we will give several examples of those $k$-algebras. 

\par \vspace{2mm}
Let us explain the organization of our paper. 
In Section 2, we recall several definitions:
Buchsbaum property, initial degree, regularity, and $a$-invariant. 
We also recall two important theorems for Buchsbaum 
homogeneous $k$-algebras: Hoa--Miyazaki theorem, 
Eisenbud--Goto theorem. 
These theorems will play key roles in our proof of the main theorem.   
\par \vspace{2mm}
In section 3, we give a lower bound on multiplicities for 
Buchsbaum homogeneous $k$-algebras in terms of 
its dimension, codimension and initial degree. 
Namely, we prove: 

\par \vspace{2mm}\par \noindent 
{\bf Theorem \ref{Ineq} (Lower bound on multiplicities of Buchsbaum $k$-algebras).}
Let $A$ be a homogeneous Buchsbaum $k$-algebra with $d=\dim A \ge 1$, 
$c=\codim A \ge 1$ and $q=\indeg A \ge 2$.  
Let $\fraM$ denote the unique graded maximal ideal of $A$. 
Let $Q$ be a parameter ideal of $A$ which is generated 
by a linear system of parameters.  
Then 
\begin{enumerate}
\item $\Sigma(Q) \subseteq \fraM^{q-1}+Q$. 
\item The following inequality holds$:$
\[
 e(A) \ge \bbinom{c+q-2}{q-2} + 
\sum_{i=1}^{d-1} \bbinom{d-1}{i-1} l_A(H_{\fraM}^i(A)).
\]
\end{enumerate}
Moreover, the equality holds in $(1)$ $($or $(2)$$)$
 if and only if $[\Sigma(Q)]_n = A_n$ for all $n \ge q-1$. 

\par \vspace{2mm}
The ring $A$ is said to be a 
\textit{Buchsbaum $k$-algebra with minimal multiplicity of degree $q$}
if the equality holds. 
Note that this notion coincides with having minimal multiplicity defined 
by the first author (resp. by Terai and the second author) when 
$q=2$ (resp. for Stanley--Reisner rings). 
Moreover, those algebras are \textit{not} Cohen--Macaulay except 
polynomial rings. In fact, we know that 
$e(A) \ge \bbinom{c+q-1}{q-1}$ if $A$ is Cohen--Macaulay; 
see, e.g., Proposition \ref{CM-LB}.
 
\par \vspace{2mm}
In Section 4, we prove the following theorem, 
which gives a characterization of Buchsbaum $k$-algebras with 
minimal multiplicity of degree $q$. 

\par \vspace{2mm} \par \noindent 
{\bf Theorem \ref{Main}(Characterization of Buchsbaum $k$-algebras with minimal 
multiplicity of degree $q$).} 
Under the same notation as in Theorem \ref{Ineq}, 
the following conditions are equivalent$:$ 
\begin{enumerate}
 \item $A$ has minimal multiplicity of degree $q$. 
 \item $a(A) = q-d-2$. 
 \item $H_{\fraM}^i(A) = [H_{\fraM}^i(A)]_{q-1-i}$ for all $i <d$ and 
 $[H_{\fraM}^d(A)]_n = 0$ for all $n \ge q-d-1$. 
 \item $A$ has $q$-linear resolution with 
\[
 \sum_{i=0}^{d-1} \bbinom{d}{i}  l_A(H_{\fraM}^i(A)) 
= \bbinom{\reg A+c-1}{c-1}. 
\]
 \item $\Sigma (Q) = \fraM^{q-1}+Q$, 
that is, $[\Sigma(Q)]_n =  A_n$ $(n \ge q-1)$. 
\end{enumerate}
When this is the case, we have 
\begin{enumerate}
\item[(a)] $\reg A =q-1$. 
\item[(b)] $\Soc(H_{\fraM}^d(A)) = [H_{\fraM}^d(A)]_{q-d-2}$. 
\end{enumerate}

\par \vspace{2mm}
In particular, the above theorem yields that any Buchsbaum $k$-algebra 
with minimal multiplicity of degree $q$ has $q$-linear resolution. 
Furthermore, we recognize that such a algebra attains an upper bound 
on the inequality in Hoa-Miyazaki theorem, Kamoi--Vogel theorem, 
respectively; see Remark \ref{Kamoi}. 

\par 
We also give some examples of Buchsbaum $k$-algebras $A$ having 
minimal multiplicity with $\depth A=0$. 
For example, let $S$ be a polynomial ring over a field $k$ 
with $\fraN$ the graded maximal ideal of $S$. 
If $S/J$ is a homogeneous Cohen--Macaulay $k$-algebra 
with $(q-1)$-linear resolution, then $A=S/\fraN J$ is 
a Buchsbaum $k$-algebra with minimal multiplicity of degree $q$ 
with $\depth A=0$.

\par \vspace{2mm}
In Section 5, we give several examples of Buchsbaum $k$-algebras 
with minimal multiplicity of degree $q$. 
We first show that 
any non-Cohen--Macaulay, Buchsbaum $k$-algebra with $e(A) \le 2$ 
has minimal multiplicity of degree at most $3$. 
Note that a Buchsbaum $k$-algebra with $e(A) = 2$ 
does \textit{not} necessarily have minimal multiplicity 
in the sense of Goto \cite{Go3}. 

\par 
On the other hand, one can find many examples of 
reduced Buchsbaum $k$-algebras with 
minimal multiplicity of higher degree among Stanley--Reisner rings;
see also \cite{TeYo1}. 
On the other hand, we have no examples of those algebras which are 
integral domains over an algebraically closed field.

\medskip
\section{Preliminaries}

\par 
We first briefly recall the definition of Buchsbaum rings, 
which was introduced by St\"uckrad and Vogel.  
Let $(A,\frm)$ be a Noetherian local ring with the maximal ideal $\frm$, and 
$M$ a finitely generated $A$-module of $\dim M = s$. 

\begin{defn}[{\rm see, e.g., \cite{SV}}]
The module $M$ is \textit{Buchsbaum} if it satisfies one of 
the following equivalent conditions: 
\begin{enumerate}
 \item The difference $l_A(M/QM) - e_Q(M)$ is independent of 
the choice of a parameter ideal $Q$ for $M$, which we denote by $I_A(M)$ 
and we call the Buchsbaum invariant of $M$. 
 \item Any system of parameters for $M$ forms a weak $M$-sequence: 
For any system $a_1,\ldots, a_s$ of parameters for $M$ and 
for any integer $1 \le i \le s$, one has the equality 
\[
 (a_1,\ldots,a_s)M : a_i = (a_1,\ldots,a_{i-1})M : \frm.
\]
 \item Any system of parameters for $M$ forms a d-sequence for $M$ 
in the sense of Huneke, 
that is, the equality 
\[
  (a_1,\ldots,a_{i-1})M : a_ia_j = (a_1,\ldots,a_{i-1})M : a_j 
\]
holds for all integers $1 \le i \le j \le s$. 
\end{enumerate}
\par 
We say that $A$ is a \textit{Buchsbaum ring} if $A$ is a Buchsbaum module over itself. 
The module $M$ is \textit{Cohen--Macaulay} if  it is Buchsbaum with $I_A(M) =0$.  
\end{defn} 
\par 
Now let $M$ be a Buchsbaum $A$-module of $\dim M =s$. 
Then the local cohomology module 
$H_{\frm}^i(M)$ is a finite-dimensional vector space over $A/\frm$ 
for all $i\ne s$. 
In particular, $\frm  H_{\frm}^i(M) =0$ for  all $i\ne s$ and 
$M_{\frp}$ is a Cohen--Macaulay $A_{\frp}$-module for all prime 
$\frp \in \Supp M \setminus \{\frm\}$. 
Furthermore, the following formula is known for Buchsbaum invariant: 
\[
  I_A(M) = \sum_{i=0}^{s-1} \bbinom{s-1}{i} l_A(H_{\frm}^i(M)). 
\]
 
Moreover, if $a_1\ldots,a_r \in \frm$ is part of a system of parameters for $M$,  
then the quotient module $\overline{M} =M/(a_1,\ldots,a_r)M$ is also a 
Buchsbaum $A$-module and $I_A(M) = I_A(\overline{M})$. 

\par \vspace{2mm}
The main target in this paper are homogeneous Buchsbaum $k$-algebras with linear 
resolutions.   
In the following, let us recall several fundamental definitions
 of homogeneous $k$-algebras.  
Let $k$ be an infinite field, and 
let $S=k[X_1,\ldots,X_v]$ be a homogeneous polynomial ring 
with $n$ variables over $k$ with $\deg X_i =1$. 
Let $A=S/I$ be a $d$-dimensional homogeneous $k$-algebra, 
where $I$ is a graded ideal such that $0 \ne I \subseteq (X_1,\ldots,X_v)^2S$.
Take a graded minimal free resolution over $S$:
\[
 0 \to \bigoplus_{j \in \Z} S(-j)^{\beta_{p,j}(A)}
\stackrel{\varphi_{p}}{\longrightarrow} \cdots
\stackrel{\varphi_{2}}{\longrightarrow}
\bigoplus_{j \in \Z} S(-j)^{\beta_{1,j}(A)}
\stackrel{\varphi_{1}}{\longrightarrow}
S \to A \to 0,
\] 
where $\beta_{i,j}(A)$ denotes the graded Betti number of $A$
and $S(-j)$ denotes the graded free $S$-module with 
grading given by $[S(-j)]_n = S_{n-j}$.   
Then the \textit{initial degree} of $A$ is defined by 
\[ 
 \indeg A = \min\{j \in \mathbb{Z} \,:\, \beta_{1,j}(A) \ne 0\}. 
\] 
Let $\indeg I$ be the smallest degree of minimal system of generators of $I$. 
Note that $\indeg S/I = \indeg I$ provided that 
$0 \ne I \subseteq (X_1,\ldots,X_v)^2S$. 
Similarly, the \textit{$($Castelnuovo--Mumford$)$ regularity} 
is defined by 
\[
\reg A = \max\{j-i \in \mathbb{Z} \,:\, \beta_{i,j}(A) \ne 0\}. \\
\]
Then $\reg A \ge \indeg A -1$, and $A$ has \textit{$q$-linear resolution} if 
the equality holds and $\indeg A =q$. 
Note that we consider $\indeg A=1$ and $\reg A = 0$ provided 
that $A(=S)$ is a polynomial ring. 
\par 
Put $\fraM = (X_1,\ldots,X_v)A$, the unique homogeneous maximal ideal of $A$. 
Let $H_{\fraM}^i(A)$ denote the $i$ th local cohomology module of $A$. 
This module has a natural graded structure.  
We recall the \textit{$a$-invariant} of $A$: 
\[
 a(A) = \max\{n \in \bbZ \,:\, [H_{\fraM}^d(A)]_n \ne 0 \}. 
\]
One can compute the regularity of $A$ as follows (see \cite{EiGo}):
\[
 \reg A = \min\{n \in \bbZ \,:\, [H_{\fraM}^i(A)]_j = 0 \quad 
 \text{for all $i +j > n$}\}.
\]
In particular, $a(A) + d \le \reg A$ holds, and 
the equality holds if $A$ is Cohen--Macaulay. 

\par 
The homogeneous $k$-algebra  $A$ is \textit{Buchsbaum}
 if the local ring $A_\fraM$ is Buchsbaum. 
Furthermore, we use $I(A)$ in replace of $I(A_{\fraM})$. 
As for regularity for Buchsbaum $k$-algebras, 
the following theorem is known$:$

\begin{thm}[{\rm Hoa--Miyazaki \cite{HoMi}}] \label{HoaMiyazaki}
Suppose that  $A$ is a homogeneous Buchsbaum $k$-algebra with $d=\dim A$.   
Then $\reg A \le a(A) +d +1$. 
\par
In particular, $\reg A = a(A)+d$, or $a(A)+d+1$. 
\end{thm}

\begin{cor} \label{a-invLB}
Under the same notation as in Theorem $\ref{HoaMiyazaki}$, we have 
$a(A) \ge q-d-2$.  
\end{cor}

\begin{proof}
By definition and Theorem \ref{HoaMiyazaki}, we have 
\[
 q-1 = \indeg A-1 \le \reg A \le a(A) + d+1. 
\]
Thus we get the required inequality. 
\end{proof}

\par
Now let $A=S/I$ be a homogeneous Buchsbaum $k$-algebra, and 
let $a_1,\ldots, a_d \in A_1$ be a linear system of parameters of $A$ and fix it. 
For an integer $\ell \ge 1$ we set
\[
 Q := (a_1,\ldots,a_d)A,\qquad 
\Sigma(\underline{a^{\ell}}) := 
 \sum_{i=1}^{d} (a_1^{\ell},\ldots,\widehat{a_i^{\ell}},\ldots,a_d^{\ell}): 
a_i^{\ell} +(a_1^{\ell},\ldots,a_d^{\ell}). 
\]
In particular, we write $\Sigma(Q) = \Sigma(\underline{a})$ for simplicity. 
Then $Q$ is a minimal reduction of $\fraM$, that is, $\fraM^{r+1} = Q\fraM^r$ for 
some integer $r \ge 0$. 
The following result plays a key role in this paper. 

\begin{thm}[{\rm see \cite[Theorem 4.1]{Go3}}]  \label{Mul-Sigma}
Suppose that  $(A,\frm)$ is a Buchsbaum local ring or a homogeneous $k$-algebra,
and that $Q$ is a minimal reduction of $\frm$. 
Then we have
\[
 e(A) = e_Q(A) 
= l_A(A/\Sigma(Q)) + \sum_{i=1}^{d-1} \bbinom{d-1}{i-1}  l_A(H_{\frm}^i(A)).
\]
\end{thm}

\par
Suppose that $A=S/I$ has $q$-linear resolution. 
Then one can take a minimal system of generators of $I$ 
whose degree is equal to $q$ and thus  
$I \subseteq \fraN^q$, where $\fraN = (X_1,\ldots,X_v)S$. 
Moreover, $A/Q$ also has $q$-linear resolution, that is, 
$A/Q \cong k[Y_1,\ldots,Y_c]/(Y_1,\ldots,Y_c)^q$, 
where $c=\codim A$. 
In the case of Buchsbaum rings, the following criterion for having $q$-linear resolution 
is known:

\begin{thm}[Eisenbud--Goto {\rm \cite[Corollary 1.15]{EiGo}}] \label{Eisen}
Suppose that  $A$ is a homogeneous Buchsbaum $k$-algebra with $\dim A =d$. 
Then the following conditions are equivalent$:$
\begin{enumerate}
 \item $A$ has $q$-linear resolution. 
 \item $H_{\fraM}^i(A) = [H_{\fraM}^i(A)]_{q-1-i}$ for all $\,i< d$ and 
$[H_{\fraM}^d(A)]_{n} = 0$ for  all $\,n \ge q-d$. 
 \item $\fraM^q = Q\fraM^{q-1}$. 
\end{enumerate}
\end{thm}

\medskip 
\section{Lower bound on multiplicities}

\par 
In \cite{Go3}, the first author proved an inequality
\begin{equation} \label{eq:min}
 e(A) \ge 1+ \sum_{i=1}^{d-1} \bbinom{d-1}{i-1} l_A(H_{\frm}^i(A))
\end{equation}
for any $d$-dimensional Buchsbaum local ring $(A,\frm)$ and called it a 
\textit{Buchsbaum local ring with minimal multiplicity} if the equality holds.   
The main purpose of this section is to prove an improved version of the above inequality 
for homogeneous Buchsbaum $k$-algebras in terms of $d=\dim A$, $c=\codim A$ and 
$q=\indeg A$ and to define the notion of  
\textit{Buchsbaum homogeneous $k$-algebras with minimal multiplicity  of degree $q$}. 

\par \vspace{2mm}
Throughout this section, we use the following notation, unless otherwise specified:
Let $S=k[X_1,\ldots,X_v]$ is a polynomial ring and $I$ is an ideal of $S$ with 
$q:=\indeg I \ge 2$, and 
let $A=S/I$ be a homogeneous $k$-algebra 
with $d:=\dim A \ge 1$, $c:=\codim A \ge 1$ and $\indeg A =q$, 
Let $\fraM$ (resp. $\fraN$) denote the unique graded maximal ideal of $A$ (resp. $S$). 
Also, let $Q$ be a parameter ideal of $A$ which is generated 
by a linear system of parameters $a_1,\ldots,a_d$. 
\par
Before stating our result, let us recall a lower bound on multiplicities 
for Cohen--Macaulay $k$-algebras. 
Under the above notation, we have 
\[
 l_A(A/Q) \ge \dim_k k[Y_1,\ldots,Y_c]/(Y_1,\ldots,Y_c)^q= \bbinom{c+q-1}{q-1}
\]
since $A/Q$ is isomorphic to $k[Y_1,\ldots,Y_c]/L$, where $L$ is an ideal 
such that  $L \subseteq (Y_1,\ldots,Y_c)^q$. 
Further, the equality holds if and only if $a(A/Q) =q-1$.  
From this one can easily obtain the following (see also \cite{EiGo}). 

\begin{prop} \label{CM-LB}
Under the above notation, 
we also suppose that $A$ is Cohen--Macaulay. 
Then 
\begin{enumerate}
 \item $e(A) \ge \bbinom{c+q-1}{q-1}$.  
 \item $a(A) = \reg A -d \ge q-d-1$.  
 \item The following conditions are equivalent$:$ 
\begin{enumerate}
 \item $e(A) = \bbinom{c+q-1}{q-1}$. 
 \item $a(A) = q-d-1$.
 \item $A$ has $q$-linear resolution.
\end{enumerate}
When this is the case, $A/Q \cong k[Y_1,\ldots, Y_c]/(Y_1,\ldots,Y_c)^q$,  
where $k[Y_1,\ldots, Y_c]$ is a polynomial ring in $c$ variables over $k$.  
\end{enumerate} 
\end{prop}

\par \vspace{2mm}
The following theorem is the first main result in this paper, 
which generalizes the inequality in Eq.(\ref{eq:min}) and 
in Proposition \ref{CM-LB}. 

\begin{thm} \label{Ineq}  
Under the above notation, we also suppose that $A$ is Buchsbaum. 
Then 
\begin{enumerate}
\item $\Sigma(Q) \subseteq \fraM^{q-1}+Q$. 
\item The following inequality holds$:$
\[
 e(A) \ge \bbinom{c+q-2}{q-2} + \sum_{i=1}^{d-1} \bbinom{d-1}{i-1} 
l_A(H_{\fraM}^i(A)).
\]
\end{enumerate}
Moreover, the equality holds in $(1)$ $($or $(2)$$)$
 if and only if $[\Sigma(Q)]_n = A_n$ for all $n \ge q-1$. 
\end{thm}

\begin{proof}
(1) It suffices to show that $[\Sigma(Q)]_n = [Q]_n$ for all $n \le q-2$. 
Fix $i$ with $1 \le i \le d$. 
Set $B = A/(a_1,\ldots,\widehat{a_i},\ldots,a_d)A$. 
Then $B$ is a $1$-dimensional Buchsbaum $k$-algebra, and 
we can write as in the form $B=k[Y_1,\ldots,Y_{c+1}]/I_B$ and 
$H_{\fraM}^0(B) = U/I_B$ for some indeterminates $Y_1,\ldots,Y_{c+1}$ and 
some ideals $I_B \subseteq U$ of $k[Y_1,\ldots,Y_{c+1}]$. 
Then $\indeg U \ge q-1$ since 
$\fraM U \subseteq I_B$ and $\indeg I_B \ge \indeg I =q$. 
Therefore if $n \le q-2$ then 
\[
\left[
\dfrac{(a_1,\ldots,\widehat{a_i},\ldots,a_d):a_i}{(a_1,\ldots,\widehat{a_i},\ldots,a_d)} \right]_n 
= [H_{\fraM}^{0}(B)]_n =0. 
\]
Namely, $[(a_1,\ldots,\widehat{a_i},\ldots,a_d):a_i]_n = 
[(a_1,\ldots,\widehat{a_i},\ldots,a_d)]_n \subseteq [Q]_n$ for every $n$. 
This implies that $[\Sigma(Q)]_n \subseteq [Q]_n$ for all $n \le q-2$. 
\par  
(2) By Theorem \ref{Mul-Sigma}, it is enough to show 
$l_A(A/\Sigma(Q)) \ge \bbinom{c+q-2}{q-2}$. 
By (1), we have
\[
 l_A(A/\Sigma(Q)) \ge \sum_{n=0}^{q-2} \dim_k [A/\Sigma(Q)]_n = 
\sum_{n=0}^{q-2} \dim_k [A/Q]_n.  
\] 
Since $\indeg A/Q \ge \indeg A=q$, we get 
\[
\sum_{n=0}^{q-2} \dim_k [A/Q]_n = \sum_{n=0}^{q-2} \bbinom{c+n-1}{n} = \bbinom{c+q-2}{q-2}. 
\]
It is easy to see the last assertion. 
\end{proof}

\vspace{1mm} 
\begin{remark} \label{indeg-remark}
Under the same notation as in Theorem \ref{Ineq}, 
if, in addition, $q > q' \ge 2$, then  
\[
 e(A) > \bbinom{c+q'-2}{q'-2} + \sum_{i=1}^{d-1} \bbinom{d-1}{i-1} 
 l_A(H_{\fraM}^i(A)).
\] 
Hence one can replace $\indeg A=q$ with $\indeg A \ge q$ in the above theorem. 
\end{remark}

\par \vspace{2mm}
By virtue of the above theorem, we can define the notion of 
Buchsbaum homogeneous $k$-algebras with minimal multiplicity, 
which generalizes the notion defined by Goto in \cite{Go3}.

\begin{defn}[\textbf{Minimal multiplicity}] 
\label{Minimal-defn}
Under the same notation as in Theorem \ref{Ineq},  
the ring $A$ is said to be a 
Buchsbaum homogeneous $k$-algebras with \textit{minimal multiplicity of degree $q$} 
if the equality holds:
\[
 e(A) = \bbinom{c+q-2}{q-2} 
+ \sum_{i=1}^{d-1} \bbinom{d-1}{i-1} \, l_A(H_{\fraM}^i(A)). 
\]
We regard a polynomial ring as a Buchsbaum (Cohen--Macaulay) 
$k$-algebra with minimal multiplicity of degree $1$. 
\end{defn}

\par \vspace{2mm}
For Buchsbaum local rings,  
we can prove a similar inequality as that in the above theorem. 
Namely, we have: Let $(A,\frm,k)$ be a  $d$-dimensional 
Buchsbaum local ring with $c = \codim A$ and suppose that 
the $\frm$-adic completion $\widehat{A}=R/I$, 
where $(R,\frn)$ is a complete regular local ring and 
$I$ is an ideal of $R$ such that $I \subseteq \frn^q$. 
Also, let $Q$ be a minimal reduction of $\frm$. 
Then the same assertions in Theorem \ref{Ineq} hold;
see the proof for \cite{GoYo}. 
Note that this yields \cite[Theorem 2.2]{HeIk} and \cite[Theorem 4.1]{Go3}
when $q=2$. 
\par
However we do \textit{not} know whether the associated graded ring 
$\gr_{\frm}(A) := \oplus_{n \ge 0} \frm^n/\frm^{n+1}$ of 
a Buchsbaum local ring with minimal multiplicity of degree $q$ 
also has the same property (especially Buchsbaum property) 
except the case of $q=2$ yet. 
Namely, we do not have a complete answer to the following question below; 
See also \cite{Go2}. 
So we do \textit{not} introduce the notion of 
\textit{Buchsbaum local ring with minimal multiplicity of degree $q$} here.

\begin{quest} \label{Graded}
Let $(A,\frm,k)$ be a $d$-dimensional Buchsbaum local ring with 
$\codim A =c$. 
Suppose that $\widehat{A} = R/I$ such that $I \subseteq \frn^q$, 
where $(R,\frn)$ is a complete regular local ring and $I$ is an ideal of $R$. 
\par 
If $\Sigma(Q) = \frm^{q-1}+Q$ $($or, $\frm^q=Q\frm^{q-1}$$)$,  
then is the associated graded ring $\gr_{\frm}(A)$ of $A$ Buchsbaum?
\end{quest}

\par
The inequality in Theorem \ref{Ineq} does \textit{not} necessarily hold 
in general for non-Buchsbaum rings.

\begin{exam} \label{nonBbm}
Let $k$ be any field, and let $q$, $d \ge 2$ be integers. 
\par \noindent 
Put $A=k[X_0,X_1,\ldots,X_d]/(X_0(X_0,\ldots,X_d)^{q-1})$. 
Then 
\begin{enumerate}
 \item $A$ is a homogeneous $k$-algebra with $\dim A =d$, $\codim A=1$ 
and $\indeg A =q$ with the unique maximal ideal $\frm=(x_0,\ldots,x_d)A$. 
 \item $H_{\frm}^0(A) = x_0A$ and $A/H_{\frm}^0(A) \cong k[X_1,\ldots,X_d]$. 
In particular, $e(A)=1$. 
 \item $A$ is Buchsbaum if and only if $q=2$. 
 \item $e(A) \ge \bbinom{c+q-2}{q-2} + 
\sum_{i=1}^{d-1} \bbinom{d-1}{i-1} \,l_A(H_{\frm}^i(A))=q-1$ 
if and only if $q=2$. 
\end{enumerate}
\end{exam}

\medskip
\section{Characterization}

In this section, we give several equivalent conditions of 
Buchsbaum homogeneous $k$-algebras with minimal multiplicity 
of degree $q$. 
In particular, we show that those rings have $q$-linear resolutions and 
satisfy several boundary conditions; e.g., Hoa--Miyazaki theorem and 
Kamoi--Vogel theorem.   

\par
In this section, we use the same notation as in the previous section, unless 
othwise specified. 
The following theorem is the second main result in this paper. 

\begin{thm} \label{Main}
Let $A$ be a homogeneous Buchsbaum $k$-algebra with $d=\dim A \ge 1$, 
$c=\codim A \ge 1$ and $q=\indeg A \ge 2$.  
Let $Q$ be a parameter ideal of $A$ which is generated 
by a linear system of parameters.  
Then 
the following conditions are equivalent$:$ 
\begin{enumerate}
 \item $A$ has minimal multiplicity of degree $q$, that is, 
\[
 e(A)  =  \bbinom{c+q-2}{q-2} + \sum_{i=1}^{d-1} \bbinom{d-1}{i-1}  
l_A(H_{\fraM}^i(A)).
\] 
 \item $a(A) = q-d-2$. 
 \item $H_{\fraM}^i(A) = [H_{\fraM}^i(A)]_{q-1-i}$ for all $i <d$ and 
 $[H_{\fraM}^d(A)]_n = 0$ for all $n \ge q-d-1$. 
 \item $A$ has $q$-linear resolution with 
\[
 \sum_{i=0}^{d-1} \bbinom{d}{i}  l_A(H_{\fraM}^i(A)) 
= \bbinom{\reg A+c-1}{c-1}. 
\]
 \item $\Sigma (Q) = \fraM^{q-1}+Q$, 
that is, $[\Sigma(Q)]_n =  A_n$ for all $n \ge q-1$. 
\end{enumerate}
When this is the case, we have 
\begin{enumerate}
\item[(a)] $\reg A =q-1$. 
\item[(b)] $\Soc(H_{\fraM}^d(A)) = [H_{\fraM}^d(A)]_{q-d-2}$. 
\end{enumerate}
\end{thm}

\par \vspace{2mm}
In the following, let us prove Theorem \ref{Main}. 
Before proving the theorem, let us prove the following lemma. 

\begin{lemma} \label{Sigma}
Suppose that $A$ is a Buchsbaum $k$-algebra with $q$-linear resolution. 
Then $[A/\Sigma(Q)]_{q-1} \cong [H_{\fraM}^d(A)]_{q-d-1}$. 
\end{lemma}

\begin{proof}
Now consider the following direct system:
\[
 \left(A/\Sigma(\underline{a})\right)(d) \longrightarrow 
 \left(A/\Sigma(\underline{a}^2)\right)(2d) \longrightarrow 
\cdots \longrightarrow 
\left(A/\Sigma(\underline{a}^{\ell})\right)(\ell d) 
\longrightarrow 
\cdots, 
\]
where the natural map 
$\left(A/\Sigma(\underline{a}^{k})\right)(k d) 
\longrightarrow 
\left(A/\Sigma(\underline{a}^{\ell})\right)(\ell d)$ 
is given by a multiplication of $(a_1 \cdots a_d)^{\ell-k}$ for $1 \le k < \ell$. 
Then $\varphi_\ell \colon \left(A/\Sigma(\underline{a}^\ell)\right)(\ell d)
\longrightarrow H_{\fraM}^d(A)$ is injective and 
$\varinjlim  \left(A/\Sigma(\underline{a}^{\ell})\right)(\ell d)
= H_{\fraM}^d(A)$ since $A$ is Buchsbaum; see \cite{Go2} or \cite[Lemma 1.3]{Ya}. 
The lemma follows from the following claim. 

\par \vspace{1mm}  \par \noindent 
{\bf Claim.} $A_{q-1+(\ell -1)d} \subseteq 
(a_1\cdots a_d)^{(\ell-1)}\fraM^{q-1} + \Sigma(\underline{a}^{\ell})$ for every 
$\ell \ge 1$. 
\par \vspace{2mm} 
In fact, the claim follows from 
\[
 \fraM^{q-1+(\ell-1)d} = 
Q^{(\ell-1)d}\fraM^{q-1} \subseteq 
(a_1^{\ell},\ldots,a_d^{\ell})\fraM^{q-1+(\ell-1)d-\ell} 
+ (a_1\cdots a_d)^{\ell-1}\fraM^{q-1}. 
\]
\end{proof}

\begin{proof}[Proof of Theorem $\ref{Main}$]
We have already seen $(1) \Longleftrightarrow (5)$ in 
the proof of Theorem \ref{Ineq}. 
So it suffices to show  
$(5) \Longrightarrow (4) \Longrightarrow (3)  \Longrightarrow 
(2) \Longrightarrow (5)$. 
Now we put $h^i(A) = l_A(H_{\fraM}^i(A))$ for all $i <d$. 

\par \vspace{1mm} 
$(5) \Longrightarrow (4):$
The assumption implies that $\fraM^q \subseteq \fraM \Sigma(Q) \subseteq Q$. 
This means that $\fraM^q = Q\cap \fraM^{q} = Q\fraM^{q-1}$ 
since $A$ is homogeneous. 
Thus $A$ has $q$-linear resolution by Eisenbud--Goto theorem (Theorem \ref{Eisen}).
In particular, $\reg A = q-1$, and 
we get
\begin{equation} \label{EQlinear}
  e(A) = e_Q(A) = l_A(A/Q) - I(A) 
= \bbinom{c+q-1}{q-1} - \sum_{i=0}^{d-1} \binom{d-1}{i} h^i(A).
\end{equation}
\par
On the  other hand, by $(5) \Longrightarrow (1)$ we have 
\begin{equation} \label{EQminimal}
  e(A) = \bbinom{c+q-2}{q-2} + \sum_{i=1}^{d-1} \binom{d-1}{i-1} h^i(A).
\end{equation}
Comparing two equalities, we get 
\begin{equation} \label{EQgoto}
  \sum_{i=0}^{d-1} \bbinom{d}{i} h^i(A) = \bbinom{c+q-2}{c-1} = 
\bbinom{\reg A +c-1}{c-1}.
\end{equation}
\par 
$(4) \Longrightarrow (3):$ 
By a similar argument as above, we obtain that 
$A$ has minimal multiplicity of degree $q$. 
In particular, $[\Sigma(Q)]_{q-1} = A_{q-1}$. 
This implies that $[H_{\fraM}^d(A)]_{q-1-d} =0$ by Lemma \ref{Sigma}
since $A$ has $q$-linear resolution. 
On the other hand, 
$H_{\fraM}^i(A) = [H_{\fraM}^i(A)]_{q-1-i}$ for all $i <d$ and 
$[H_{\frm}^d(A)]_{n}=0$ for all $n \ge q-d$ 
by Eisenbud--Goto theorem. 
Hence we get (3).   

\par \vspace{1mm}  
$(3) \Longrightarrow (2):$
By the assumption we have $a(A) \le q-d-2$. 
However the converse is always true; see Corollary \ref{a-invLB}. 

\par \vspace{1mm}  
$(2) \Longrightarrow (5):$
Assume that $a(A) = q-d-2$. Then for all $n \ge q-1$, we have 
\[
 \left[A/\Sigma(Q)\right]_n \stackrel{\varphi_1}{\hookrightarrow} 
[H_{\fraM}^d(A)]_{n-d} = 0.
\]
Namely, $[\Sigma(Q)]_n = A_n$ for all $n \ge q-1$, as required. 

\par \vspace{2mm}
We now prove that $\Soc(H_{\fraM}^d(A)) = [H_{\fraM}^d(A)]_{q-d-2}$
if $a(A) = q-d-2$. 
It is enough to show that $[\Soc(H_{\fraM}^d(A))]_{n} = 0$ for 
all $n \le q-d-3$. 
Consider the following direct system described as above$:$
\[
 [A/\Sigma(Q)]_{n+d} \stackrel{a_1\cdots a_d}{\longrightarrow}
 [A/\Sigma(\underline{a}^2)]_{n+2d} \stackrel{a_1\cdots a_d}{\longrightarrow}
\cdots \longrightarrow [H_{\fraM}^d(A)]_{n}.
\]
By \cite[Proposition 3.8]{Ya}, we have 
\[
 \Soc(H_{\fraM}^d(A)) \subseteq \Hom_A(A/Q,H_{\fraM}^d(A))
 = \varphi_2\left(
\frac{\sum_{i=1}^d a_1\cdots \widehat{a_i} \cdots a_d \Sigma(Q)
+\Sigma(\underline{a}^2)}{\Sigma(\underline{a}^2)}
\right).
\]
Since $\Sigma(Q) = \fraM^{q-1}+Q$, if $j \le q+d-3$, then 
\begin{eqnarray*}
\left[\sum_{i=1}^d a_1\cdots \widehat{a_i}\cdots a_d \Sigma(Q) 
+ \Sigma(\underline{a}^2)\right]_j 
&=& \sum_{i=1}^d a_1\cdots \widehat{a_i}\cdots a_d [Q]_{j-d+1}
+ [\Sigma(\underline{a}^2)]_j  \\
&\subseteq& a_1\cdots a_d A_{j-d} + [\Sigma(\underline{a}^2)]_j.
\end{eqnarray*}
Hence $[\Soc(H_{\fraM}^d(A))]_n \subseteq [\Image \varphi_1]_n$, 
that is,  $[\Soc(H_{\fraM}^d(A))]_n = 
\left[\frac{\Sigma(Q):\fraM}{\Sigma(Q)}\right]_{n+d}$ for all $n \le q-d-3$. 
It is enough to show the following claim. 

\par \vspace{1mm} \par \noindent 
{\bf Claim.} $[\Sigma(Q):\fraM]_j \subseteq [\Sigma(Q)]_j$ for all $j \le q-3$. 
\par \vspace{1mm}
Fix $j \le q-3$. 
For any element $\xi \in [\Sigma(Q):\fraM]_j$, 
$\fraM(\fraM \xi) \subseteq \fraM \Sigma(Q) \subseteq Q$. 
Then for any element $a \in A_1$, 
$a\xi \in [\fraM \xi]_{j+1} \subseteq [Q:\fraM]_{j+1} = [Q]_{j+1}$
because $\Soc(A/Q)$ is concentrated in degree $q-1$. 
Hence $\fraM \xi \subseteq Q$. 
That is, $\xi \in [Q:\fraM]_{j} = [Q]_j \subseteq [\Sigma(Q)]_j$, as required. 
\par \vspace{1mm}
We have finished the proof of Theorem \ref{Main}. 
\end{proof}

\par \vspace{1mm}
The first author \cite{Go4} determined the Hilbert series of Buchsbaum 
$k$-algebras with $q$-linear resolution. 
Thus we can also describe that of Buchsbaum 
$k$-algebras with minimal multiplicity. 

\begin{cor} \label{Hilbert}
Suppose that $A$ is a Buchsbaum homogeneous $k$-algebras 
with minimal multiplicity of degree $q$. 
Put $d=\dim A \ge 2$ and $c=\codim A \ge 1$. 
Then 
The Hilbert series $F(A,t)$ of $A$ is given as follows$:$  
\begin{eqnarray*}
F(A,t) &= &  
\dfrac{1}{(1-t)^d} \sum_{i=0}^{q-1} \bbinom{c+i-1}{i} t^i  \cr
 & & 
 + \dfrac{1}{(1-t)^d} \sum_{i=q}^{q+d-1} (-1)^{i+q-1} 
\left\{\sum_{j=0}^{q+d-1-i} \bbinom{d}{i+j-q+1} h^j(A)\right\} t^i. 
\end{eqnarray*}
\end{cor}

\par \vspace{1mm}
Roughly speaking, a Buchsbaum homogeneous $k$-algebra 
with minimal multiplicity is \textit{not} Cohen--Macaulay. 

\begin{cor} \label{non-CM}
Suppose that $A$ is a Buchsbaum homogeneous $k$-algebras 
with minimal multiplicity of degree $q$. 
Then $A$ is Cohen--Macaulay if and only if it is isomorphic 
to a polynomial ring and $q=1$. 
\end{cor}

\begin{proof}
It follows from Proposition \ref{CM-LB}(2) and Theorem \ref{Main}(2).  
\end{proof}

\vspace{1mm}
In general, $A/aA$ does \textit{not} necessarily have 
minimal multiplicity even if so is $A$ and $a \in A_1$ is a non-zero-divisor. 

\begin{cor} \label{non-NZD}
Suppoose that $A$ is a $d$-dimensional Buchsbaum homogeneous $k$-algebras 
with minimal multiplicity of degree $q$. 
Let $a \in A_1$ be a non-zero-divisor of $A$ and set $\overline{A} = A/aA$. 
Then 
\begin{enumerate}
\item $\overline{A}$ is a Buchsbaum $k$-algebra with $q$-linear resolution. 
\item $\overline{A}$ has minimal multiplicity 
if and only if $H_{\fraM}^{d-1}(A)=0$. 
\end{enumerate}
\end{cor}

\begin{proof}
In the following, we put $h^i(A) = l_A(H_{\fraM}^i(A))$ for each $i \le d-1$. 
\par 
(1) The assumption implies that $A$ has $q$-linear resolution 
by Theorem \ref{Main} $(1)\Longrightarrow (4)$.   
Hence $\overline{A}$ also has $q$-linear resolution.
\par \vspace{1mm} \par 
(2) We know that $h^i(\overline{A})= h^i(A) + h^{i+1}(A)$ for each $i \le d-2$
since $A$ is Buchsbaum and  $a \in \fraM$ is a non-zero-divisor. 
Hence 
\[
\sum_{i=1}^{d-2} \bbinom{d-2}{i-1} h^i(\overline{A})
= \sum_{i=1}^{d-1} \bbinom{d-1}{i-1} h^i(A) - h^{d-1}(A). 
\]
The assertion easily follows from this. 
\end{proof}

\par \vspace{1mm}
Let us consider the case where $A/H_{\fraM}^0(A)$ is Cohen--Macaulay. 
Then one can obtain many examples of Buchsbaum homogeneous $k$-algebras with 
minimal multiplicity of degree $q$. 

\begin{cor} \label{WCM}
Under the same notation as in Theorem $\ref{Main}$, 
if, in addition, $A/H_{\fraM}^0(A)$ is Cohen--Macaulay, 
then the following conditions are equivalent$:$
\begin{enumerate}
 \item $A$ has minimal multiplicity of degree $q$, that is, 
$e(A) = \bbinom{c+q-2}{q-2}$.  
 \item $A$ has $q$-linear resolution 
and $l_A(H_{\fraM}^0(A)) = \bbinom{c+q-2}{q-1}$. 
 \item $H_{\fraM}^0(A) = [H_{\fraM}^0(A)]_{q-1}$, 
$H_{\fraM}^i(A) = 0$ for all $1 \le i \le d-1$ and 
$[H_{\fraM}^d(A)]_n = 0$ for all $n \ge q-d-1$. 
 
 \item $A/H_{\fraM}^0(A)$ has $(q-1)$-linear resolution. 
\end{enumerate}
\end{cor}

\begin{proof}
It follows from Theorem \ref{Main} and Proposition \ref{CM-LB}. 
\end{proof}

\begin{remark} \label{saibu}
The condition (4) does \textit{not} 
imply that $A$ has minimal multiplicity
(we need to assume $\indeg A=q$);
see also Example \ref{depth-zero}.  
\end{remark}

\begin{exam} 
 If $S/U$ is a Cohen--Macaulay $k$-algebra with $(q-1)$-linear resolution, 
then $A=S/\fraN U$ is a Buchsbaum $k$-algebra with minimal multiplicity of
degree $q$ with $l_A(H_{\fraM}^0(A)) = \mu_S(U)$, where $\mu_S(U)$ denotes 
the minimal number of system of generators of $U$. 
For example, 
\begin{enumerate}
\item $A=S/f \fraN$ is a Buchsbaum $k$-algebra with 
minimal multiplicity of degree $q$ provided that 
$f \in \fraN^{q-1} \setminus \fraN^q$. 
 \item For given integers $c \ge 1$, $q \ge 2$, if we put 
\[
 A=k[X_1,\ldots,X_c,Y_1,\ldots,Y_d]/
(X_1,\ldots,X_c,Y_1,\ldots,Y_d)(X_1,\ldots,X_c)^{q-1}, 
\]  
then $A$ is a Buchsbaum $k$-algebra with minimal multiplicity of degree $q$ 
such that $\dim A=d$, $\codim A = c$ and $A/H_{\fraM}^0(A)$ is Cohen--Macaulay. 
\end{enumerate}
\end{exam}

\vspace{2mm}
\begin{remark} \label{Kamoi}
In \cite{KaVo}, Kamoi and Vogel proved that the following inequality$:$
\[
 \sum_{i=0}^{d-1} \bbinom{d}{i}\,h^i(A) \le \binom{\reg A+c-1}{c-1}
\]
for any homogeneous Buchsbaum $k$-algebra $A$.  
Our theorem \ref{Main}(3) shows that any Buchsbaum ring with minimal multiplicity 
attains the upper bound in this inequality. 
We do \textit{not} know whether the above equality implies 
that $A$ has linear resolution
(and thus it has minimal multiplicity in our sense).  
\par 
Furthermore, those rings satisfy  \lq\lq $\reg A = a(A) + d+1$'' because 
$\reg A = q-1$ and $a(A) = q-d-2$. See Hoa--Miyazaki theorem. 
\end{remark}
 
\begin{remark} \label{2-linear}
For a Buchsbaum $k$-algebra $A$,  
it has maximal embedding dimension if and only if 
it has $2$-linear resolution. 
\par 
It is known that any Buchsbaum ring with minimal multiplicity 
has $2$-linear resolution; see \cite{Go3}. 
Our theorem \ref{Main} $(1) \Longleftrightarrow (4)$ generalizes this fact. 
\end{remark}

\medskip
\section{Examples}

\par 
In this section, we first show that 
all non-Cohen--Macaulay, Buchsbaum $k$-algebras with $e(A) \le 2$ 
have minimal multiplicities of degree at most $3$. 
Next, we give several examples of Buchsbaum homogeneous $k$-algebras 
with minimal multiplicity of higher degree using 
the theory of Stanley--Reisner rings; see \cite{TeYo1}. 
\par
Throughout this section, let $k$ be a field and 
$S=k[X_1,\ldots, X_v]$ a polynomial ring over $k$ 
with the unique graded maximal ideal $\fraN$. 

\subsection{Buchsbaum rings with small multiplicities}

\vspace{2mm}
In \cite{Go1}, the first author classified non-Cohen--Macaulay, 
Buchsbaum rings $(A,\frm)$ with $e(A)=2$ and proved that 
such a ring has minimal multiplicity in the sense of Goto \cite{Go3} 
such that $h^i(A):l_A(H_{\frm}^i(A))=0$ for all $i \ne 1,d$ 
and $h^1(A) = 1$ when $d \ge 2$ 
and $\depth A > 0$; see \cite[Theorem 1.1]{Go1} 
and also \cite[Section 4]{Go3} for details. 
But, when $\depth A=0$, there exist a non-Cohen--Macaulay, Buchsbaum ring
with $e(A)=2$ which does not have minimal multiplicity
in the sense of Goto \cite{Go3}. 
Now let us show that such a ring has 
minimal multiplicity in our sense. 
Before doing it, we prove the following lemma. 

\begin{lemma} \label{IndegMult}
Suppose that $A=S/I$ is a homogeneous non-Cohen--Macaulay, 
Buchsbaum $k$-algebra. 
Put $d =\dim A\ge 2$, $c=\codim A \ge 1$, $q=\indeg A \ge 2$ and $e=e(A)$. 
Then 
\begin{enumerate}
 \item $q \le e+1$. 
 \item The following conditions are equivalent$:$
\begin{enumerate}
 \item $q=e+1$. 
 \item There exists an element $f \in \fraN^{q-1} \setminus \fraN^{q}$ 
such that $I=f\fraN$. 
 \item $A$ has minimal multiplicity of degree $e+1$. 
\end{enumerate}
\end{enumerate}
When this is the case, $c=1$ and $h^i(A)=0$ for all $i=1,\ldots,d-1$. 
\end{lemma}

\begin{proof}
(1) Suppose that $q \ge e+2$. 
By Theorem \ref{Ineq} we have 
\[
 e\ge \bbinom{c+q-2}{q-2} \ge \bbinom{c+e}{e} \ge \bbinom{e+1}{e}=e+1. 
\]
This is a contradiction. 
Hence $q \le e+1$.  
\par 
(2) It is enough to show $(a) \Longrightarrow (b)$. 
We may assume that $e \ge 2$. 
Suppose (a).
By Theorem \ref{Ineq} again, we have 
\[
 e  \ge  \bbinom{c+e-1}{e-1} + \sum_{i=1}^{d-1} \bbinom{d-1}{i-1} h^i(A) 
\ge \bbinom{e}{e-1} + \sum_{i=1}^{d-1} \bbinom{d-1}{i-1} h^i(A) 
\ge e. 
\]
This implies that $c=1$, $h^i(A)=0$ for all $i\ne 0,d$ and that 
$A$ has minimal multiplicity of degree $q = e+1$. 
Take an ideal $U$ such that $H_{\fraM}^0(A) = U/I$. 
Then $\indeg U \ge 2$ since $\fraN U \subseteq I$. 
Hence $\codim A/H_{\fraM}^0(A) =1$ and thus 
$A/H_{\fraM}^0(A)$ is a hypersurface. 
So one can take an element $f \in \fraN$ such that $U=fS$. 
Since $f\fraN \subseteq I \subseteq fS$ and $e(A/H_{\fraM}^0(A))=e=q-1$, 
we obtain (b). 
\end{proof}

\begin{remark} \label{Known-lemma}
Most part of the above lemma is probably known. 
In fact, Hoa--Miyazaki theorem implies that 
$q-1 \le \reg A \le a(A)+d+1 \le e(A)$ since $A$ is Buchsbaum. 
\end{remark}

\begin{prop} \label{MinMulti2}
Suppose $d \ge 2$. 
Then any homogeneous non-Cohen--Macaulay, Buchsbaum $k$-algebra $A$ with 
$e(A) \le 2$ has  minimal multiplicity whose degree is at most $3$. 
To be precise, we have$:$ 
\begin{enumerate}
 \item Suppose that $e(A)=1$.  
Then $A$ has minimal multiplicity of degree $2$ with $\depth A=0$. 
Moreover, $A/H_{\fraM}^0(A)$ is isomorphic to a polynomial ring over $k$. 
 \item Suppose that $e(A) =2$. 
Then $A$ has minimal multiplicity of degree $q$ for some integer $q \ge 2$, and 
one of the following three cases occurs$:$ 
\begin{enumerate}
\item $q=2$ and $\depth A > 0$. Then $c=d$ and $h^1(A)=1$. 
\item $q=2$ and $\depth A = 0$. Then $h^1(A)=1$. 
\item $q=3$ and $\depth A = 0$. Then $c=1$ and $A/H_{\fraM}^0(A)$ is Cohen--Macaulay. 
\end{enumerate}
\end{enumerate}
\end{prop}

\begin{proof}
When $\depth A > 0$, the assertion follows from \cite{Go1,Go3}. 
So we may suppose that $\depth A = 0$ and that $q \ge 2$. 
\par 
(1) The assertion is clear.  See also Theorem \ref{Ineq}. 
\par 
(2) By Lemma \ref{IndegMult}, we have  $q \le 3$. 
\vspace{1mm}\par \noindent 
{\bf Case 1.} when $h^1(A) \ge 1$: 
\par
We obtain that $q=2$ by Lemma \ref{IndegMult}. 
Then $A$ has minimal multiplicity and $h^1(A)=1$ since $e(A) = 2 = 1+h^1(A) \ge 2$. 
\vspace{1mm}\par \noindent 
{\bf Case 2.} when $h^1(A) =0$: 
\par
Set $A=S/I$, $\fraM = \fraN/I$, $W:=H_{\fraM}^0(A):=U/I$ and $B:= A/W = S/U$. 
Then $B$ is a hypersurface because it is Cohen--Macaulay and $e(B)=2$. 
Thus there exists $f \in \fraN^2 \setminus \fraN^3$ such that $U=fS$. 
Since $f \fraN \subseteq I \subseteq fS$, we have $I=f\fraN$. 
In particular, $q=3$. 
Thus $A$ has minimal multiplicity of degree $3$ and $c=1$ by Lemma \ref{IndegMult}.
\end{proof}

\begin{exam}[See also \cite{Go1}] \label{EXmulti2}
Let $X_1,\ldots,X_d,Y_1,\ldots,Y_d,Y,Z_1,\ldots,Z_c$ be indeterminates
 over a field $k$. 
\begin{enumerate}
 \item Put $A=k[X_1,\ldots,X_d,Y]/(X_1Y,\ldots,X_dY,Y^2)$. 
Then $A$ is a Buchsbaum $k$-algebra with $\dim A = d$, $\depth A =0$ and $e(A)=1$. 
Also, it has minimal multiplicity of degree $2$. 
 \item Put $A=k[X_1,\ldots,X_d,Y_1,\ldots,Y_d]/(X_1,\ldots,X_d)\cap(Y_1,\ldots,Y_d)$. 
Then $A$ is a Buchsbaum $k$-algebra with $\dim A = d$, $\depth A = 1$ and $e(A)=2$. 
Also, it has minimal multiplicity of degree $2$. 
 \item Put $S=k[X_1,\ldots,X_d,Y_1,\ldots,Y_d,Z_1,\ldots,Z_c]$.  
Also, if we put $I=(X_1,\ldots,X_d)\cap(Y_1,\ldots,Y_d)+(Z_1,\ldots,Z_c)\fraN$ 
and $A=S/I$, 
then $A$ is a Buchsbaum $k$-algebra with $\dim A = d$, $\depth A =0$, $\codim A=c$ and  
$e(A)=2$. 
Also, it has minimal multiplicity of degree $2$. 
 \item Put $A = k[X_1,\ldots,X_d,Y]/(X_1Y^2,\ldots,X_dY^2)$. 
Then $A$ is a Buchsbaum $k$-algebra with $\dim A = d$, $\depth A =0$ and 
$e(A)=2$. 
Also, it has minimal multiplicity of degree $3$. 
\end{enumerate}
\end{exam}

\par  \vspace{2mm}
By a similar argument as in the proof of Proposition \ref{MinMulti2}
one can show that any non-Cohen--Macaulay Buchsbaum $k$-algebra $A$ with 
$e(A)=3$ has minimal multiplicity in our sense whenever $q \ge 3$; see 
also Examples \ref{Hyp}, \ref{depth-zero}. 
However, this is \textit{not} true in general when $q=2$; 
see Example \ref{depth-zero}(2). 

\begin{exam}[See also Lemma \ref{IndegMult}] \label{Hyp}
Let $U$ be a graded ideal of $S$ such that $B=S/U$ is Cohen--Macaulay 
with $e(B) =3$. 
If we put $A=S/\fraN U$, then $A$ is Buchsbaum such that 
$A/H_{\fraM}^0(A) \cong B$ and $e(A)=e(B)=3$. 
Then either one of the following two cases occurs:
\begin{enumerate}
 \item $B$ has maximal embedding dimension and $A$ has minimal multiplicity of 
degree $3$. 
 \item $B$ is a hypersurface and $A$ has minimal multiplicity of 
degree $4$. 
\end{enumerate}
\end{exam}

\par
There exists a Buchsbaum $k$-algebra $A$ having 
minimal multiplicity with $\depth A=0$ such that $A/H_{\fraM}^0(A)$ is 
not Cohen--Macaulay; see the next example (3). 

\vspace{1mm}
\begin{exam} \label{depth-zero} 
Let $S = k[X,Y,Z,W]$ be a polynomial ring. 
\begin{enumerate}
\item If we set $I = (X,Y,Z,W)(XW-YZ, Y^2-XZ, Z^2-YW)$ and $A = S/I$,   
then $A/H_{\fraM}^0(A) \cong k[s^3,s^2t,st^2,t^3]$ is Cohen--Macaulay with  
$2$-linear resolution. 
Thus $A$ is Buchsbaum with minimal multiplicity of degree $3$. 
Also, we have$:$ 
\[
 H_{\fraM}^0(A) = k(-2), \quad 
 H_{\fraM}^1(A) = 0\;\; \text{and}\;\; e(A) =3.
\]
\item If we set $I=(X,Y,Z,W)(YZ,Z^2) + (Y^2-XZ)$ and $A = S/I$, 
then $A/H_{\fraM}^0(A) 
\cong k[X,Y,Z,W]/(Y^2-XZ,YZ,Z^2) \cong \gr_{\frm}(k[[t^3,t^4,t^5]])[W]$ 
is Cohen--Macaulay with $e(B)=3$. 
Thus $A$ is Buchsbaum with $\indeg A (=\indeg B)=2$ and $e(A)=3$. 
But $A$ does \textit{not} have minimal multiplicity and 
\[
 H_{\fraM}^0(A) = k(-2)^{\oplus 2}, \quad 
 H_{\fraM}^1(A) = 0\;\;\text{and}\;\; e(A) =3.
\]
\item If we set $I = ((X,Y,Z,W)(XW-YZ), Z^3-YW^2,Y^3-X^2Z,XZ^2-Y^2W)$ and $A=S/I$, 
then 
$A$ is Buchsbaum with minimal multiplicity of degree $3$
since \[
H_{\fraM}^0(A) = k(-2), \quad H_{\fraM}^1(A) = k(-1) \;\; \text{and}\;\; 
e(A) =4.
\]
But $A/H_{\fraM}^0(A) \cong k[s^4,s^3t,st^3,t^4]$ is \textit{not} Cohen--Macaulay. 
\end{enumerate}
\end{exam}

\vspace{2mm}
\subsection{The case of Stanley--Reisner rings} 
One can find abundant examples of Buchsbaum reduced $k$-algebras 
with minimal multiplicity 
in the class of Stanley--Reisner rings. 
In fact, in \cite{TeYo1} 
Terai and the second author \cite{TeYo1} gave another definition of 
Buchsbaum $k$-algebras with minimal multiplicity for Stanley--Reisner rings
and have characterized those rings. 
This paper was inspired their work. 

\par \vspace{2mm}
Now let us briefly recall the theory of Stanley--Reisner rings. 
Let $\Delta$ be a \textit{simplicial complex}  
on the vertex set $V=[v]:=\{1,\ldots,v\}$, that is, 
$\Delta$ is a collection of subsets of $V$ such that 
(a) $F \subseteq G,\, G \in \Delta \Longrightarrow F \in \Delta$ 
and (b) $\{i\} \in \Delta$ for $i \in V$. 
The dimension of $F \in \Delta$ 
($F$ is said to be a \textit{face} of $\Delta$) is defined by 
$\#(F)-1$. 
Set $\dim \Delta = \max\{\dim F\,:\, F \in \Delta\}$. 
A complex is called \textit{pure} if all facets (maximal faces) have the same 
dimension. 
Put $\link_{\Delta} G = \{F \in \Delta \,:\, F \cup G 
\in \Delta,\, F \cap G \notin \Delta \}$, the \textit{link} of a face $G$ in $\Delta$.  
\par 
Set $I_{\Delta} = (X_{i_1}\cdots X_{i_p} \,:\, 1 \le i_1 <\cdots < i_p \le v,\, 
\{i_1,\ldots,i_p\} \notin \Delta)S$. 
Then $k[\Delta] = S/I_{\Delta}$ is called the \textit{Stanley--Reisner ring}
 of $\Delta$ over $k$. 
The ring $k[\Delta]$ is a homogeneous reduced $k$-algebra with 
$d:=\dim k[\Delta] = \dim \Delta+1$. 
Moreover, $q=\indeg k[\Delta] \le d+1$ always holds, and 
$\Delta$ is $(d-1)$-skeleton of $(v-1)$-simplices of $2^{[v]}$ whenever $q=d+1$.  
Put $f_i = \#\{F \in \Delta\,:\, \dim F = i\}$ for all $i \in \bbZ$. 
Then the Hilbert series of $k[\Delta]$ can be written as in the following form:
\[
 F(k[\Delta],t) = \sum_{i=-1}^{d-1} \dfrac{f_i t^{i+1}}{(1-t)^{i+1}} 
= \dfrac{h_0+h_1t+\cdots + h_dt^d}{(1-t)^d}. 
\]
In particular, $e(k[\Delta])$ is  
equal to $f_{d-1}(\Delta)$, the number of facets $F$ with $\dim F = d-1$.   

\par \vspace{1mm}
In the following, let $\Delta$ be a simplicial complex on $V=[v]$ 
with $\dim \Delta = d-1$. 
For Buchsbaumness of Stanley--Reisner rings, 
the following fact is well-known. 
Note that $k[\Delta]$ is Buchsbaum if and only if 
it has (F.L.C.), that is, every local cohomology modules 
$H_{\fraM}^i(k[\Delta])$ has finite length for $i \ne d$. 

\begin{fact} \label{BbmSR}
The following conditions are equivalent$:$ 
\begin{enumerate}
 \item $k[\Delta]$ is Buchsbaum. 
 \item $\Delta$ is pure and 
$k[\link_{\Delta} \{i\}]$ is Cohen--Macaulay $($of dimension $d-1$$)$ for 
each $i \in V$. 
 \item $H_{\fraM}^i(A) = [H_{\fraM}^i(A)]_0 (\cong \widetilde{H}_{i-1}(\Delta;k))$ 
for all $i < d$.   
\end{enumerate}
\end{fact} 

\par \vspace{2mm}
In the class of Buchsbaum Stanley--Reisner rings, there exists a criterion for 
$k[\Delta]$ to have linear resolution in terms of $h$-vectors 
as follows:

\begin{thm}[\textrm{\cite[Theorem 1.3]{TeYo1}}] \label{Linear-SR}  
Suppose that $A = k[\Delta]$ is 
Buchsbaum and $2 \le q \le d$. 
Put $h = \dim_k H^{q-1}_{\fraM}(A)$. 
Then the following conditions are equivalent$:$
\begin{enumerate}
\item $A$ has $q$-linear resolution. 
\item The $h$-vector 
$h(\Delta) = (h_0,h_1,\ldots,h_{q-1},h_q,h_{q+1},\ldots,h_d)$ of $\Delta$ is 
\[
\left(1,\,c,\,\cdots, \sbinom{c+q-2}{q-1},\,- \sbinom{d}{q}h,
  \,\sbinom{d}{q+1}h,\cdots,(-1)^{d-q+1} \sbinom{d}{d}h\right)
\]
 \item The following equalities hold$:$
\[
e(A)= \bbinom{c+q-1}{q-1} - h \bbinom{d-1}{q-1}
 \quad \text{and}
 \quad I(A) = h \bbinom{d-1}{q-1}. 
\]
\end{enumerate}
When this is the case, $H_{\fraM}^i(A) = 0$ $(i \ne q-1,d)$ and 
the following inequalities hold$:$
\[
0 \;\le\; h \;\le\;
 \frac{c(c+1)\cdots (c+q-2)}{d(d-1)\cdots (d-q+2)}=:h_{c,d,q}.
\]
\end{thm}

\par \vspace{2mm}
From this point of view, one can regard a 
\textit{Buchsbaum simplicial complex with minimal multiplicity}
 as a Buchsbaum complex with linear resolution 
and \lq\lq \textit{maximal homology}''.

\begin{thm}[\textrm{\cite[Theorem 2.3]{TeYo1}}] \label{BbmSRmm}
Suppose that $A = k[\Delta]$ is Buchsbaum. 
Then 
\[
\text{\small $e(A) \ge \frac{c+d}{d} \bbinom{c+q-2}{q-2}$}
\]
holds, and the following statements are equivalent$:$
 \begin{enumerate}
  \item $\Delta$ is a Buchsbaum complex with 
minimal multiplicity of degree $q$, 
that is, $e(A) = \frac{c+d}{d} \bbinom{c+q-2}{q-2}$.  
  \item $A$ has minimal multiplicity of degree $q$ in our sense. 
  \item $A$ has $q$-linear resolution with $l_A(H_{\fraM}^{q-1}(A)) = h_{c,d,q}$. 
 \item $k[\link_{\Delta}\{i\}]$ has $(q-1)$-linear resolution for all $i \in V$. 
 \item The $h$-vector of $A$ can be written as in the form of 
\textrm{Theorem $\ref{Linear-SR}(2)$} and $h=h_{c,d,q}$. 
 \item $k[\Delta^{*}]$ has pure and almost linear resolution and 
$a(k[\Delta^{*}])=0$, where 
$\Delta^{*} = \{F \in \Delta\,:\, V \setminus F \notin \Delta\}$ denotes 
the Alexander dual complex of $\Delta:$ 
\[
 0 \to S(-(c+d))^{\beta_{q}^{*}} \to  
S(-(c+q-2))^{\beta_{q-1}^{*}} \to \cdots \to S(-c)^{\beta_{1}^{*}} 
 \to S \to k[\Delta^{*}] \to 0,    
\]
where $S = k[X_1,\ldots,X_v]$ and $\beta_q^{*} = h_{c,d,q}$.  
\end{enumerate}
\end{thm}

\par \vspace{2mm}
Let us pick up some examples of Buchsbaum complexes with minimal multiplicity
from \cite{TeYo1}. 

\begin{exam} \label{simplices}
$\Delta$ is a Buchsbaum complex with 
minimal multiplicity of degree $2$
if and only if it is a finite disjoint union of $(d-1)$-simplices. 
\par 
When this is the case, $k[\Delta]$ is isomorphic to 
\[
 k[X_{ij}: 1 \le i \le d, 1\le j \le e]/
(X_{pr}X_{qs}: 1 \le p,q \le d,\, 1\le r \ne s \le e), 
\]
where $e=e(k[\Delta])$, the number of connected components of $\Delta$. 
\end{exam}

\par \vspace{1mm}
Now let us recall a cyclic polytope. 
The algebraic curve $M \subseteq \bbR^f$ parametrically by 
$x(t) = (t,t^2,\ldots,t^f)$
is called the moment curve. 
Let $n \ge f+1$ be an integer. 
A \textit{cyclic polytope} with $n$ vertices, 
denoted by $C(n,f)$, is the convex hull of any $n$ 
distinct points on $M$. 
Note that $C(n,f)$ is a simplicial $f$-polytope. 
Next example provides us $d$-dimensional Buchsbaum 
reduced $k$-algebras of degree $q$ for given integers $2 \le q \le d$. 

\begin{exam} \label{cyclic}
Let $q$, $d$ be given integers with $2 \le q \le d$. 
Put $n = 2d-q+2$ and $f = 2(d-q+1)$. 
Let $\Delta$ be the Alexander dual of 
the boundary complex $\Gamma$ of 
a cyclic polytope $C(n,f)$. 
Then since $k[\Gamma]$ satisfies the condition $(6)$ in Theorem \ref{BbmSRmm}, 
we obtain that $k[\Delta]$ is a $d$-dimensional 
Buchsbaum Stanley--Reisner ring with minimal multiplicity 
of degree $q$ with $h=h_{c,d,q}=1$. 
\end{exam}

\vspace{1mm}
For a given integer $n \ge 5$ such that $n \equiv 0,2 \pmod{3}$, 
there exists a $2$-dimensional Buchsbaum complex on $[n]$ 
with minimal multiplicity of degree $3$; see \cite{Ha, TeYo1, TeYo2}.  

\begin{exam} \label{hanano}
Let $\Delta$ be a simplicial complex on $V=[5]$ spanned by 
$\{1,2,3\}$, $\{1,3,4\}$, $\{1,4,5\}$, $\{2,3,5\}$ and $\{2,4,5\}$
as follows$:$

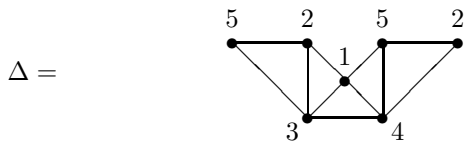
\setlength{\unitlength}{1mm}
\begin{figure}[ht] \label{Fig:SRring}
\begin{picture}(100,22)
\put(0,10){$\Delta=$}
\put(40,5){\line(1,0){10}} 
\put(40,5){\line(0,1){10}} 
\put(40,5){\line(-1,1){10}} 
\put(40,5){\line(1,1){10}} 
\put(50,5){\line(0,1){10}} 
\put(50,5){\line(1,1){10}} 
\put(40,15){\line(-1,0){10}} 
\put(40,15){\line(1,-1){10}} 
\put(50,15){\line(1,0){10}} 
\put(39,4){$\bullet$}
\put(49,4){$\bullet$}
\put(39,14){$\bullet$}
\put(49,14){$\bullet$}
\put(29,14){$\bullet$}
\put(59,14){$\bullet$}
\put(44,9){$\bullet$}
\put(44,12){$1$}
\put(39,17){$2$}
\put(49,17){$5$}
\put(29,17){$5$}
\put(59,17){$2$}
\put(37,2){$3$}
\put(51,2){$4$}
\end{picture}
\end{figure}
\par \vspace{-9mm} \par \noindent
Then 
\begin{eqnarray*}
 k[\Delta]
&=& k[X_1,X_2,X_3,X_4,X_5]/
(X_1X_2X_4, X_1X_2X_5, X_1X_3X_5, X_2X_3X_4, X_3X_4X_5) \\
& = & k[X_1,X_2,X_3,X_4,X_5]/
(X_1,X_4) \cap  (X_4,X_5) \cap  (X_2,X_5) \cap  (X_2,X_3) \cap  (X_1,X_3)
\end{eqnarray*}
is a $3$-dimensional Buchsbaum rings with minimal multiplicity of 
degree $3$. 
\end{exam}

\par
In general, the following problem remains open when $d \ge 4$. 
When $q=d$, it suffices to show the existence of 
$(d-1)$-dimensional Buchsbaum complexes with minimal multiplicity of
degree $q$ with $\codim k[\Delta]=c$; see \cite{TeYo2}.

\begin{prob} \label{Existence}
Let $c,d,q,h$ be integers with $c \ge 2$, $2 \le q \le d$ and 
$0 \le h \le h_{c,d,q}$. 
Construct $(d-1)$-dimensional Buchsbaum complexes $\Delta$
with $q$-linear resolution such that $\codim k[\Delta] = c$ and 
$\dim_k H_{\fraM}^{q-1}(k[\Delta]) = h$. 
\end{prob}

\medskip
\subsection{The case of integral domains}
There are many examples of homogeneous Buchsbaum integral domains 
having linear resolutions; see e.g., \cite{Am1, EiGo, Go4}.  
However we have no examples of Buchsbaum homogeneous integral domain over 
an algebraically closed field $k$ 
with minimal multiplicity of degree $q \ge 2$.
Indeed, our examples described before are not integral domains. 

\begin{discuss}[See \cite{EiGo}] \label{Domain-exist}
Let $d \ge 2$, $h^i$ $(1 \le i \le d-1)$, $s \ge 0$ be integers.  
Let $k$ an algebraically closed field. 
Then there exist $d$-dimensional 
Buchsbaum homogeneous $k$-algebra $A$ which are integral domains 
having $q$-linear resolutions 
with $\codim A =2$ and such that 
$\dim_k H_{\fraM}^i(A) = h^i$ for all $i=1,\ldots,d-1$, 
where 
\[
 q = \sum_{i=1}^{d-1} 
\left\{\sum_{j=1}^{d+2-i}(-1)^j \cdot (j-2) \cdot \bbinom{d+2}{i+j}\right\}
h^i + s-3. 
\]
But this number $q$ is too big 
to satisfy the condition (4) in Theorem \ref{Main}.     
\end{discuss}

\par 
The following proposition is maybe well-known for experts. 

\begin{prop} \label{two-domain}
Let $A$ be a Buchsbaum homogeneous $k$-algebra with minimal multiplicity 
of degree at most $2$. 
Suppose that $k$ is algebraically closed and that $A$ is an integral domain. 
Then $A$ is isomorphic to a polynomial ring. 
\end{prop}

\begin{proof}
By Abhyankar's result mentioned as in the introduction, 
we have $\emb(A) \le e(A) + \dim A -1$. 
On the other hand, $\emb(A) = e(A) + \dim A -1+ I(A)$ since  
$A$ has maximal embedding dimension. 
Hence $I(A) =0$, that is, $A$ is Cohen--Macaulay. 
Thus $A$ is a polynomial ring by Corollary \ref{non-CM}. 
\end{proof}

\par 
Based on the above observation we pose the following conjecture:

\begin{conj} \label{Domain}
There is no Buchsbaum homogeneous integral domain 
over an algebraically closed field with minimal multiplicity of degree $q \ge 2$. 
\end{conj}

\par
In the above conjecture we cannot remove the assumption that $k$ is algebraically closed.  

\begin{exam}[See \cite{Go1}]  \label{Exam-domains}
Let $K/k$ be an extension of fields of degree $2$, and  
set $B = K[X_1,\ldots,X_n]$, a polynomial ring.  
Then $A = \{f \in B \,:\,  f(0,\ldots,0) \in k \}$ 
is a homogeneous Buchsbaum $k$-algebra with minimal multiplicity of degree $2$
since $e(A) =2$. 
\end{exam}


\end{document}